\documentclass[reqno,11pt]{amsart}
\usepackage{amsfonts}

\usepackage{amsmath}
\usepackage{times}
\usepackage[dvips]{graphicx}

\newtheorem{theorem}{Theorem}[section]

\newtheorem{corollary}[theorem]{Corollary}

\newtheorem{lemma}[theorem]{Lemma}

\newtheorem{proposition}[theorem]{Proposition}

\newcommand{\func}[1]{\mathrm{#1}}

\begin{document}
\title{On free stochastic differential equations}
\author{V. Kargin }
\thanks{Department of Mathematics, Stanford University, Palo Alto, CA 94305,
USA. e-mail: kargin@stanford.edu}
\date{November 2010}
\maketitle

\begin{center}
\textbf{Abstract}
\end{center}

\begin{quotation}
The paper derives an equation for the Cauchy transform of the solution of a
free stochastic differential equation (SDE). This new equation is used to
solve several particular examples of free SDEs.\bigskip
\end{quotation}

\section{Introduction}

Free stochastic differential equations generalize classical stochastic
differential equations to the setting of free probability. Here is an
example of such an equation: \ 
\begin{equation*}
X_{t}=a(X_{t})dt+b^{\ast }(X_{t})(dZ_{t})b(X_{t})
\end{equation*}%
In this equation, $X_{t}$ is a self-adjoint operator, $a(X_{t})$ and $%
b(X_{t})$ are operator-valued functions of $X_{t}$, and the driving noise $%
Z_{t}$ is an operator process with free increments. That is, the increments $%
Z_{s}-Z_{t},$ $s>t,$ are assumed to be free from past realizations of $%
Z_{t}. $ The process $Z_{t}$ is usually the free Brownian process, in which
case the increments have semicircle distributions; however, other choices
are possible.

Informally, the reader may think about $X_{t}$ as very large random matrices
and $Z_{t}$ as matrices with independent Gaussian random variables as
entries. These entries follow independent Brownian motions and we are
interested in the law of the eigenvalues of $X_{t}.$ The free probability
theory is a convenient abstraction which intends to model the situation when
the size of the matrices is very large.

The study of free stochastic differential equations (``free stochastic
calculus'') is more difficult than in the classical case because of
non-commutativity of coefficients and noise. This paper contributes by
developing a new tool for the analysis of these equations.

The idea of free stochastic calculus was first suggested in \cite{speicher90}%
. It was later developed and formalized in \cite{kummerer_speicher92}, \cite%
{biane97}, \cite{biane_speicher98}, and \cite{anshelevich02}, which
introduced stochastic integration with respect to free Brownian motion as a
rigorous basis for free stochastic calculus. They also derived an analog of
the It\^{o} formula, which allows us to obtain identities like the
following:\ 
\begin{equation*}
\int_{0}^{a}[W_{t}^{2}(dW_{t})+W_{t}(dW_{t})W_{t}+(dW_{t})W_{t}^{2}]=W_{a}^{3}-2\int_{0}^{a}W_{t}dt,
\end{equation*}%
where $W_{t}$ denotes the free Brownian motion (the Wigner process). \ An
analogous formula in the classical situation is 
\begin{equation*}
\int_{0}^{a}3B_{t}^{2}(dB_{t})=B_{a}^{3}-3\int_{0}^{a}B_{t}dt,
\end{equation*}%
where $B_{t}$ is the standard Brownian motion. Note the different
coefficient before the integral on the right-hand side.

The classical It\^{o} formula is very helpful in the study of stochastic
differential equations. Unfortunately, the range of applicability of the
free It\^{o} formula is smaller. This difficulty calls for a different
method applicable to those free SDEs, which are not solvable with the It\^{o}
formula. One possibility is to seek an equation for the evolution of the
spectral distribution of the solution. Such an equation was derived by Biane
and Speicher in \cite{biane_speicher01} for the equation 
\begin{equation}
X_{t}=a(X_{t})dt+dW_{t}.  \label{SDE_biane_speicher}
\end{equation}%
They showed that the density of the spectral probability measure of $X_{t},$
which we denote $p_{t}$ and which generalizes the eigenvalue distribution of
a matrix, satisfies the free Fokker-Planck equation:%
\begin{equation}
\frac{\partial p_{t}}{\partial t}=-\frac{\partial }{\partial x}%
[p_{t}(Hp_{t}+a)].  \label{PDE_free_Fokker_Planck}
\end{equation}%
Here $H$ denotes a multiple of the Hilbert transform: 
\begin{equation}
Hu(x):=\mathrm{p.v.}\int \frac{u(y)}{x-y}dy.
\label{formula_Hilbert_transform}
\end{equation}

Still, the approach through the free Fokker-Planck equation has its own
disadvantages. First, it is applicable only to equations that have the
special form (\ref{SDE_biane_speicher}), that is, only to equations with the
constant diffusion coefficient. Second, the free Fokker-Planck equation (\ref%
{PDE_free_Fokker_Planck}) is not a bona fide partial differential equation
since it includes the Hilbert transform operator. For this reason, it is
somewhat difficult to solve this equation.

The purpose of this paper is to approach the free stochastic equations by
deriving a differential equation for the Cauchy transform of the solution.

Recall that the \emph{resolvent} of operator $X_{t}$ is defined as the
operator-valued function of a complex parameter $G_{t}\left( z\right)
:=\left( X_{t}-z\right) ^{-1}.$ The \emph{Cauchy transform} of $X_{t}$ is
defined as the expectation of the resolvent: $g_{t}(z):=E[\left(
X_{t}-z\right) ^{-1}].$ It is useful because the knowledge of the Cauchy
transform is sufficient to recover all \ properties of the spectral
probability distribution of $X_{t}$. It turns out that if $X_{t}$ solves 
\begin{equation*}
dX_{t}=a(X_{t})dt+b(X_{t})(dW_{t})c(X_{t}),
\end{equation*}%
then $g_{t}(z)$ satisfies the following equation:%
\begin{equation}
\frac{dg_{t}}{dt}%
=-E(a_{t}G_{t}^{2})+E(b_{t}c_{t}G_{t})E(b_{t}c_{t}G_{t}^{2}),
\label{equation_g_evolution_Wigner}
\end{equation}%
where we use $a_{t},$ $b_{t},$ and $c_{t}$ to denote $a(X_{t}),$ $b(X_{t}),$
and $c(X_{t}),$ respectively. This is the statement of Theorem \ref%
{theorem_g_evolution_Wigner} below.

Equation (\ref{equation_g_evolution_Wigner}) is not a usual differential
equation since it involves expectations. In general, these expectations are
difficult to compute because the coefficients $a_{t},$ $b_{t},$ and $c_{t},$
and the resolvent $G_{t}$ are not free from each other. However, if the
coefficients are polynomials, it is possible to perform further reduction to
a differential equation as we will show in Proposition \ref%
{proposition_g_evolution_PDE}.

As we just said, the knowledge of the Cauchy transform can be used to
recover the spectral probability distribution. In particular, the free
Fokker-Planck equation (\ref{PDE_free_Fokker_Planck}) can be derived from (%
\ref{equation_g_evolution_Wigner}) as will be shown in Corollary \ref%
{corollary2}.

In certain cases it is not possible to compute the Cauchy transform
explicitly, but it is possible to detect the behavior of its singularities.
This knowledge can provide us with information about the support of the
spectral distribution. In particular, it can show us how the norm of the
solution grows.

For a simple example of this approach, let us consider the well-known case
of the free Ornstein-Uhlenbeck equation: 
\begin{equation}
dX_{t}=-\theta X_{t}dt+\sigma dW_{t}.  \label{SDE_free_Ornstein_Uhlenbeck}
\end{equation}%
For this equation, it is easy to compute the Cauchy transform using equation
(\ref{equation_g_evolution_Wigner}) and recover the known result that for
positive $\theta ,$ the spectral probability distribution of $X_{t}$
converges to a stationary solution, which is a semicircle distribution
supported on the interval $[-\sigma \sqrt{2/|\theta |},\sigma \sqrt{%
2/|\theta |}].$

As a more difficult example, consider the equation 
\begin{equation*}
dX_{t}=\theta X_{t}dt+X_{t}^{1/2}( dW_{t}) X_{t}^{1/2},
\end{equation*}%
which can be thought of as a free analog of the equation for the ``geometric
Brownian motion'', $dx_{t}=\theta x_{t}dt+x_{t}dB_{t}.$

Let $X_{0}=I.$ Equation (\ref{equation_g_evolution_Wigner}) leads to the
following differential equation: 
\begin{equation*}
\frac{\partial g}{\partial t}+z(\theta -1-zg)\frac{\partial g}{\partial z}%
=-g(\theta -1-zg),
\end{equation*}%
with the initial condition $g\left( 0,z\right) =\left( 1-z\right) ^{-1}.$
The method of characteristics gives us a functional equation for the Cauchy
transform: 
\begin{equation}
z+g^{-1}=e^{(\theta -1-zg)t}.  \label{equation_geom_Brown0}
\end{equation}

\begin{figure}[tbph]
\begin{center}
\includegraphics[width=10cm]{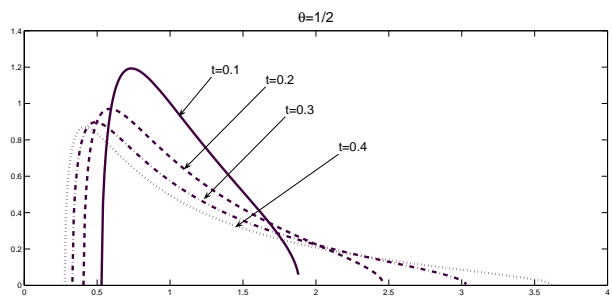}
\end{center}
\caption{$\protect\theta=1/2$}
\label{figure_gB1}
\end{figure}

\begin{figure}[tbph]
\begin{center}
\includegraphics[width=10cm]{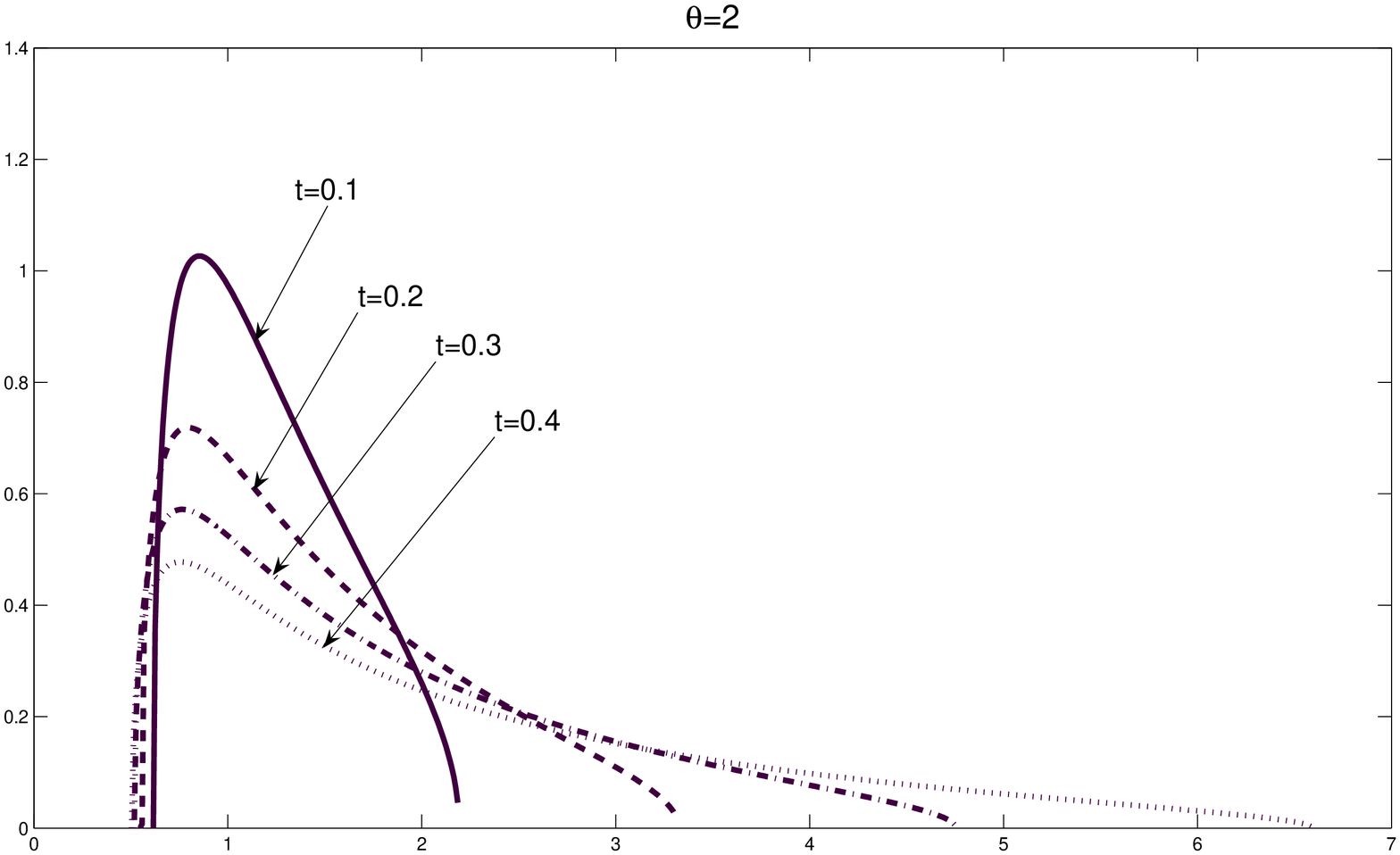}
\end{center}
\caption{$\protect\theta=2$}
\label{figure_gB2}
\end{figure}

\begin{figure}[tbph]
\begin{center}
\includegraphics[width=10cm]{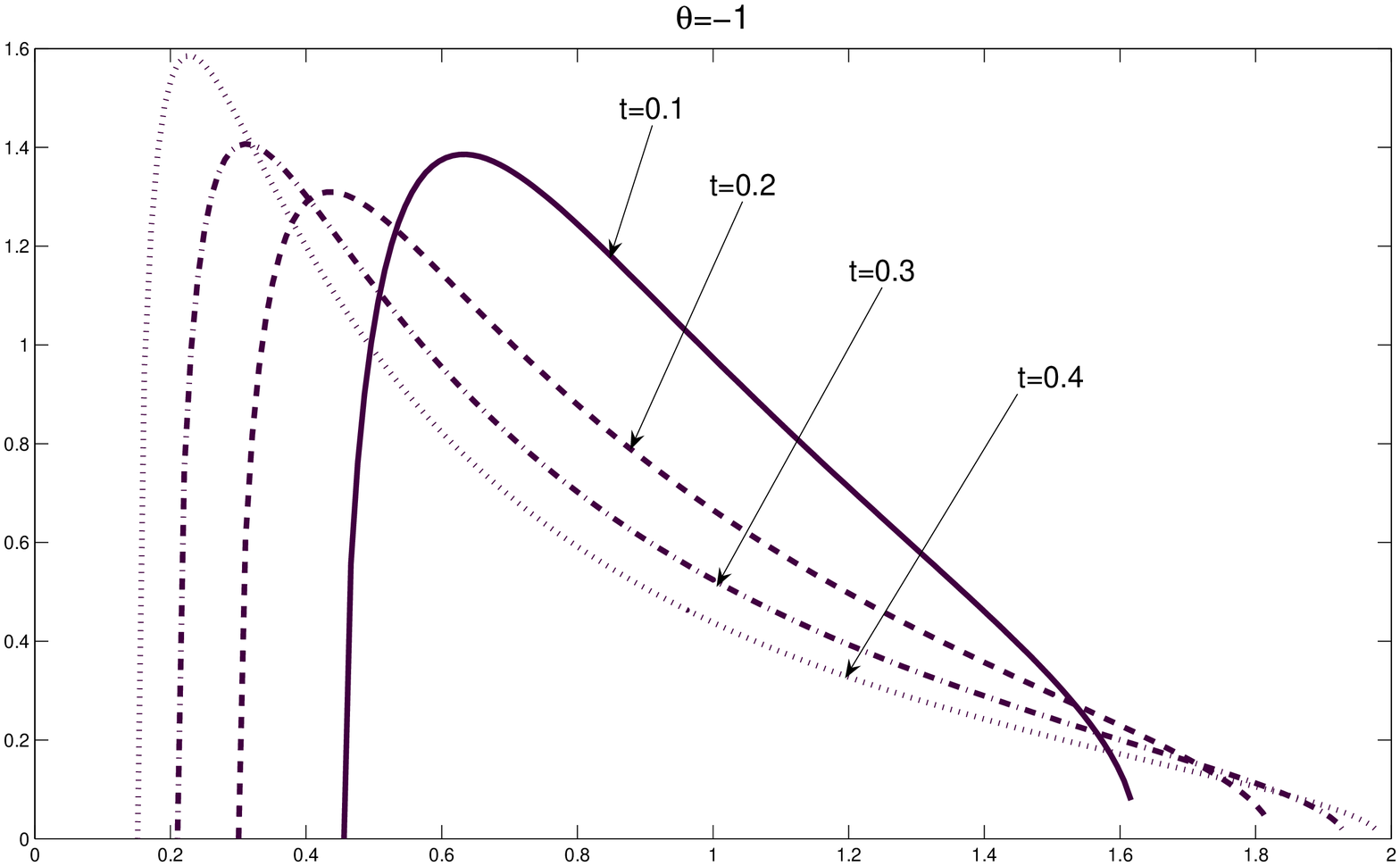}
\end{center}
\caption{$\protect\theta=-1$}
\label{figure_gB3}
\end{figure}

\begin{figure}[tbph]
\begin{center}
\includegraphics[width=10cm]{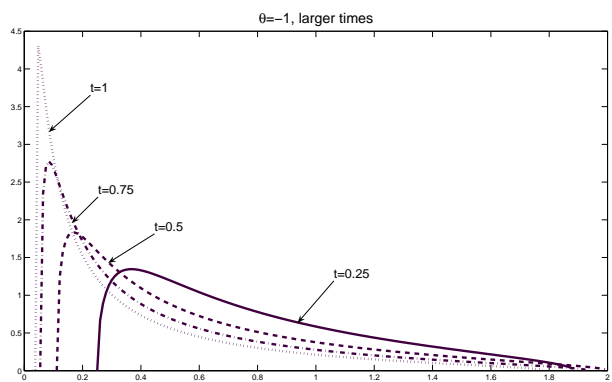}
\end{center}
\caption{$\protect\theta =-1,$ Large times}
\label{figure_gB4}
\end{figure}

While it is difficult to extract an explicit analytical formula for the
solution of this equation, we can investigate how the support of the
distribution changes with time. It turns out that for $\theta <0,$ the
support of the distribution shrinks to zero. If $\theta $ is between $0$ and 
$1,$ then the lower boundary of the support decreases to zero and the upper
boundary grows exponentially fast to infinity. If $\theta >1,$ then both the
lower and upper boundary of the support grow exponentially fast to infinity.

We can solve equation (\ref{equation_geom_Brown0}) numerically and recover
the density of the spectral distribution by using Stieltjes formula. Figures %
\ref{figure_gB1}-\ref{figure_gB4} show the evolution of the density for
various values of parameters and illustrate the complexity of the behavior
of the spectral distribution. For example, Figure \ref{figure_gB2} shows
that even if $\theta >1,$ the spectral distribution does not approach
infinity immediately. There is a transition period in which a significant
portion of the spectral distribution remains below $\lambda =1.$ Similarly,
Figure \ref{figure_gB3} shows that for $\theta <0,$ the distribution does
not collapse to zero immediately. Only when time increases, the distribution
begins the rapid approach to zero, as shown in Figure \ref{figure_gB4}.

Let us compare this result with the classical analog. By using the It\^{o}
formula, it is easy to show that the solution of the classical equation for
the geometric Brownian motion is 
\begin{equation*}
x_{t}=\exp \{(\theta -\frac{1}{2})t+B_{t}\}.
\end{equation*}%
Hence, with probability $1,$ the classical solution will decrease
exponentially to zero if $\theta <1/2,$ and will grow exponentially to
infinity if $\theta >1/2.$ However, the support of the solution distribution
is $(0,\infty )$ for all $t$. This is quite unlike the behavior of the free
SDE solution.

Note that the equation for the geometric Brownian motion can be generalized
to the free probability setting in a different way: 
\begin{equation*}
dX_{t}=\theta X_{t}dt+X_{t}dW_{t}+(dW_{t})X_{t}.
\end{equation*}%
The behavior of the solution of this equations is quite different. In
particular, the ratio of the standard deviation to the expectation is $\sqrt{%
2(e^{2t}-1)}.$ This ratio grows exponentially fast with $t,$ quite unlike
the previous example, where this ratio equals $\sqrt{t}.$ Unfortunately, the
partial differential equation associated with equation is more difficult to
solve and it is not clear whether the solution becomes unbounded in finite
time.

Finally, let us consider the following equation:

\begin{equation*}
dX_{t}=kX_{t}(dW_{t})X_{t},
\end{equation*}%
and let the initial condition be $X_{0}=aI.$ For this equation it is
possible to write an explicit formula for the spectral distribution of the
solution. An interesting feature of this equation is that the solution blows
up in finite time $\tau =(ak)^{-2},$ by which we means that the operator
norm of the solution becomes infinite as $t$ approaches $\tau $.

\begin{figure}[tbph]
\begin{center}
\includegraphics[width=10cm]{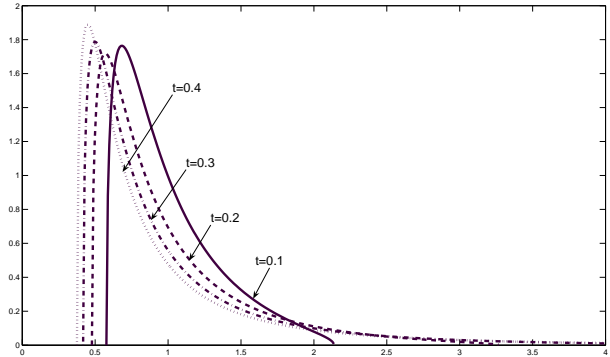}
\end{center}
\caption{Density functions for $t=0.1$, $0.2$, $0.3$, $0.4$. }
\label{figure_expl1}
\end{figure}

\begin{figure}[tbph]
\begin{center}
\includegraphics[width=10cm]{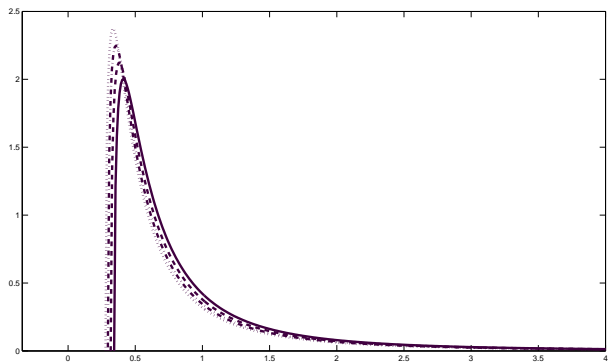}
\end{center}
\caption{Density functions for $t=0.5$, $0.6$, $0.7$, $0.8$ }
\label{figure_expl2}
\end{figure}

Another interesting feature is that as time $t$ approaches $\tau ,$ the
spectral distribution converges to a fixed distribution. If $k=a^{-1}$, then
the density of this distribution is%
\begin{equation*}
f\left( \xi \right) =\frac{\sqrt{4\xi -1}}{2\pi \xi ^{3}},
\end{equation*}%
which is supported on the interval $\left[ 1/4,\infty \right) .$ Otherwise,
it is a scaled version of this distribution. The behavior of the solution
density for various times is illustrated in Figures \ref{figure_expl1} and %
\ref{figure_expl2}.

Several specific classes of free SDE have already been investigated in the
literature. Biane and Speicher in \cite{biane_speicher01} and Gao in \cite%
{gao06} studied the free Ornstein-Uhlenbeck equation (\ref%
{SDE_free_Ornstein_Uhlenbeck}). Biane and Speicher proved that its solution
converges to a stationary process with a semicircle distribution. Gao
considered free Ornstein-Uhlenbeck processes with a free Levy driving noise,
and showed that every self-decomposable probability measure on the real line
can be realized as a distribution of such a process.

Capitaine and Donati-Matin in \cite{capitaine_donati-martin05} defined the
free Wishart process and found that it satisfies the free SDE of the form:%
\begin{equation*}
dX_{t}=\lambda dt+\sqrt{X_{t}}dZ_{t}+dZ_{t}^{\ast }\sqrt{X_{t}},
\end{equation*}%
where $Z_{t}$ is the complex Wigner process. Demni in \cite{demni08} studied
the so-called free Jacobi processes which satisfy equations similar to the
following:%
\begin{equation*}
dX_{t}=( \theta I-X) dt+\sqrt{I-X_{t}}dZ_{t}\sqrt{X_{t}}+\sqrt{X_{t}}%
dZ_{t}^{\ast }\sqrt{I-X_{t}}.
\end{equation*}

With exception of the free Ornstein-Uhlenbeck process, we study a different
set of free SDEs, and we approach these equations with a different point of
view based on the differential equations for the Cauchy transform.

For the free Ornstein-Uhlebeck process, our results agree with results in %
\cite{biane_speicher01}.

The rest of the paper is organized as follows. Section \ref%
{section_preliminaries} provides preliminary information about free
stochastic integration and It\^{o} formulas. Section \ref%
{section_dif_equations} describes main results. . In particular, Section \ref%
{subsection_existence} is devoted to a local existence and uniqueness
result. Section \ref{subsection_theory} presents general results about the
Cauchy transform of the solution and Section \ref{subsection_examples}
provides examples.

\section{Free Stochastic Integration \label{section_preliminaries}}

\subsection{The free Brownian motion}

For the basics of free probability theory we refer to \cite%
{voiculescu_dykema_nica92} and \cite{nica_speicher06}. All operators that we
consider belong to a non-commutative $W^{\ast }$-probability space $(%
\mathcal{A},E)$, that is, to a von Neumann operator algebra $\mathcal{A}$
with a faithful normal trace $E$. We denote the usual operator norm by $%
\Vert X\Vert $, and the $L^{2}$-norm by $\Vert X\Vert _{2}:=\sqrt{E(X^{\ast
}X)}$.

The \emph{spectral probability distribution} 
\index{spectral probability distribution}of a self-adjoint operator $X\in 
\mathcal{A}$ is a probability measure $\mu $ on $\mathbb{R}$ such that 
\begin{equation*}
E\left( X^{k}\right) =\int_{\mathbb{R}}x^{k}\mu \left( dx\right) .
\end{equation*}%
Its Cauchy transform is the function 
\begin{equation*}
g_{X}\left( z\right) =\int_{\mathbb{R}}%
\frac{\mu \left( dx\right) }{x-z}.
\end{equation*}%
It can be defined directly in terms of operator $X$ as the expectation of
the resolvent: $g_{X}\left( z\right) =E[G_{X}\left( z\right) ],$ where $%
G_{X}\left( z\right) :=\left( X-z\right) ^{-1}.$ The probability measure $%
\mu $ can be recovered from its Cauchy transform by the Stieltjes inversion
formula:%
\begin{equation}
\mu \left( B\right) =\frac{1}{\pi }\lim_{\varepsilon \downarrow 0}\int_{B}%
\mathrm{Im}g\left( x+i\varepsilon \right) dx,
\label{formula_Stieltjes_inversion}
\end{equation}%
provided that $B$ is Borel and $\mu \left( \partial B\right) =0.$

This fact is the starting point of our approach, since we will study the
evolution of the Cauchy transform as a tool to investigate the evolution of
the corresponding probability measure.

The most important concept in free probability theory is that of free
independence. Let $\overline{A}_{i}$ denote an arbitrary element of algebra $%
\mathcal{A}_{i}.$ The sub-algebras $\mathcal{A}_{1},\ldots ,\mathcal{A}_{n}$
of algebra $\mathcal{A}$ (and operators that generate them) are said to be 
\emph{freely independent} or \emph{free}, if the following condition holds: 
\begin{equation*}
E(\overline{A}_{i(1)}\ldots \overline{A}_{i(m)})=0,
\end{equation*}%
provided that $E(\overline{A}_{i(s)})=0$ and $i(s+1)\neq i(s)$ for every $s$%
. Two particular consequences of this definition is that (i) $E(AB)=E(A)E(B)$
if $A$ and $B$ are free, and (ii) 
\begin{equation}
E(AX_{1}AX_{2})=E(A^{2})E(X_{1})E(X_{2}),  \label{formula_4_product}
\end{equation}
if $A$ is free from $X_{1}$ and $X_{2}$ and $E(A)=0.$

The \emph{free Brownian motion}, or the \emph{Wigner process}, is a family
of operators $W_{t},$ where $t\geq 0,$ that satisfies the following
properties: (1) $W_{0}=0$; (2) the increments of $W_{t}$ are free in the
sense of Voiculescu, i.e., if $t>s,$ then $W_{t}-W_{s}$ is free from the
subalgebra $\mathcal{W}_{s}$ which is generated by all $W_{\tau }$ with $%
\tau \leq s$ , and (3) the spectral distribution of $W_{t}-W_{s}$ is
semicircle with zero expectation and variance $t-s.$

The choice of $W_{t}$ is not unique, and in the rest of the paper we assume
that a particular realization of $W_{t}$ is fixed.

\subsection{Free stochastic integral}

It\^{o}-style free stochastic integration with respect to the free Brownian
motion was defined and studied in \cite{kummerer_speicher92} and \cite%
{biane_speicher98}. Their results show that under certain assumptions on the
operator coefficients $a_{t}$ and $b_{t},$ it is possible to define the
integral 
\begin{equation*}
I=\int_{0}^{1}a_{t}(dW_{t})b_{t},
\end{equation*}%
where $W_{t}$ is the free Brownian motion.

Let us briefly recall the construction of the integral. For details of the
construction, the reader is advised to see Definition 2.2.1 and Section 3 in %
\cite{biane_speicher98} and Section 3 and Theorem 14 in \cite{anshelevich02}%
. Suppose that $a_{t}$ and $b_{t}$ are functions of $W_{\tau },$ $\tau \leq
t.$ That is, let $a_{t}$ and $b_{t}$ belong to the sub-algebra $\mathcal{W}%
_{t}$. Assume also that $\max \{\left\| a_{t}\right\| ,\left\| b_{t}\right\|
\}\leq C$ for all $t\in \left[ 0,1\right] $ and that $t\rightarrow a_{t}$
and $t\rightarrow b_{t}$ are continuous mappings in the operator norm. Let $%
t_{0},\ldots ,t_{n}$ and $\tau _{1},\ldots ,\tau _{n}$ be real numbers such
that 
\begin{equation*}
0=t_{0}\leq t_{1}\leq \ldots \leq t_{n}=1,
\end{equation*}%
and 
\begin{equation*}
0\leq \tau _{k}\leq t_{k-1}.
\end{equation*}%
We denote the set of $t_{0},\ldots ,t_{n}$ and $\tau _{1},\ldots ,\tau _{n}$
as $\Delta $. Let 
\begin{equation*}
d(\Delta )=\max_{1\leq k\leq n}(t_{k}-\tau _{k}).
\end{equation*}

Consider the sum%
\begin{equation*}
I( \Delta ) =\sum_{i=1}^{n}a_{\tau _{i}}( W_{t_{i}}-W_{t_{i-1}}) b_{\tau
_{i}}
\end{equation*}

It turns out that as $d(\Delta )\rightarrow 0,$ the sums $I(\Delta )$
converge in operator norm and the limit does not depend on the choice of $%
t_{i}$ and $\tau _{i}.$ The limit is called the free stochastic integral and
denoted as $\int_{0}^{1}a_{t}(dW_{t})b_{t}.$ An important point in the proof
of convergence is that the convergence of sums in the operator norm depends
on a free analogue of the Burkholder-Gundy martingale inequalities. 

A very useful tool in the study of stochastic integrals is the It\^{o}
formula. A free probability analogue of the It\^{o} formula was developed in %
\cite{biane_speicher98}. In terms of formal rules, it can be written as
follows: 
\begin{eqnarray}
a_{t}dt\cdot b_{t}dt &=&a_{t}dt\cdot b_{t}dW_{t}c_{t}=a_{t}dW_{t}b_{t}\cdot
c_{t}dt=0,  \notag \\
a_{t}dW_{t}b_{t}\cdot c_{t}dW_{t}d_{t} &=&E(b_{t}c_{t})a_{t}d_{t}dt.
\label{formula_ito}
\end{eqnarray}%
Note that the rule in the second line is significantly different from the
classical case. In terms of free stochastic integrals, the second rule can
be written as follows: 
\begin{eqnarray*}
\int_{0}^{1}a_{t}dW_{t}b_{t}\cdot \int_{0}^{1}c_{t}dW_{t}d_{t}
&=&\int_{0}^{1}(\int_{0}^{t}a_{\tau }dW_{\tau }b_{\tau })c_{t}dW_{t}d_{t} \\
&&+\int_{0}^{1}a_{t}dW_{t}b_{t}(\int_{0}^{t}c_{\tau }dW_{\tau }d_{\tau }) \\
&&+\int_{0}^{1}E(b_{t}c_{t})a_{t}d_{t}dt.
\end{eqnarray*}

(Compare Theorem 4.1.2. in \cite{biane_speicher98} or Proposition 8 and
Corollary 10 in \cite{anshelevich02}.) Here is an illustration (a particular
case of Proposition 4.3.2 in \cite{biane_speicher98}):

Let $W_{t}$ be the free Brownian motion and define 
\begin{equation*}
\partial (W_{t}^{n}):=W_{t}^{n-1}dW_{t}+W_{t}^{n-2}(dW_{t})W_{t}+\ldots
+(dW_{t})W_{t}^{n-1}.
\end{equation*}%
Then, 
\begin{eqnarray*}
(W_{a})^{n} &=&\int_{0}^{\alpha }d\left( W_{t}^{n}\right)  \\
&=&\int_{0}^{\alpha }\partial \left( W_{t}^{n}\right) +\int_{0}^{\alpha
}\sum_{0\leq k+l\leq n-2}W_{t}^{n-k-l-2}\left( dW_{t}\right) W_{t}^{k}\left(
dW_{t}\right) W_{t}^{l},
\end{eqnarray*}%
and it follows that%
\begin{equation}
\int_{0}^{a}\partial (W_{t}^{n})=(W_{a})^{n}-\sum_{k=0}^{\lfloor n/2\rfloor
-1}(n-2k-1)C_{k}\int_{0}^{a}W_{t}^{n-2k-2}t^{k}dt,
\label{formula_ito_formula}
\end{equation}%
where $C_{k}$ are the Catalan numbers, 
\begin{equation*}
C_{k}:=\frac{1}{k+1}\binom{2k}{k}=E[(W_{1})^{2k}].
\end{equation*}

An analogous formula for the classical It\^{o} integral with respect to the
Brownian motion $B_{t}$ is quite different:%
\begin{equation*}
\int_{0}^{a}nB_{t}^{n-1}(dB_{t})=B_{a}^{n}-\frac{n(n-1)}{2}%
\int_{0}^{a}B_{t}^{n-2}dt.
\end{equation*}%
Below, we will use the free It\^{o} formula in order to compute the moments
of variable $X_{t}$.

\section{Free Stochastic Differential Equations\label{section_dif_equations}}

\subsection{Existence and uniqueness \label{subsection_existence}}

A free stochastic differential equation (free SDE) 
\begin{equation}
dX_{t}=a(X_{t})dt+b(X_{t})(dW_{t})c(X_{t})  \label{equation_Wigner}
\end{equation}%
is a convenient shortcut notation for the following integral equation: 
\begin{equation}
X_{t}=X_{0}+\int_{0}^{t}a(X_{\tau })d\tau +\int_{0}^{t}b(X_{\tau })dW_{\tau
}c(X_{\tau }).  \label{equation_Wigner_integral}
\end{equation}%
We consider only equations with the coefficients that do not depend
explicitly on time, and we will always assume that $a(X_{t})$, $b(X_{t}),$
and $c(X_{t})$ are locally operator Lipschitz functions. (A function $f:%
\mathbb{R}\rightarrow \mathbb{C}$ is called locally operator-Lipschitz, if
it is a locally bounded, measurable function, and if for all $A>0,$ there is
a constant $K_{A}>0,$ such that 
\begin{equation*}
\Vert f(X)-f(Y)\Vert \leq K_{A}\Vert X-Y\Vert 
\end{equation*}%
for all self-adjoint operators $X$ and $Y$ with the norm less than $A.$ For
example, all polynomials are locally operator-Lipschitz.)

Equations (\ref{equation_Wigner}) and (\ref{equation_Wigner_integral}) are
particular cases of the following more general equations: 
\begin{equation}
dX_{t}=a(X_{t})dt+\sum_{i=1}^{m}b_{i}(X_{t})(dW_{t})c_{i}(X_{t})
\end{equation}%
and 
\begin{equation}
X_{t}=X_{0}+\int_{0}^{t}a(X_{\tau })d\tau
+\sum_{i=1}^{m}\int_{0}^{t}b_{i}(X_{\tau })dW_{\tau }c_{i}(X_{\tau }).
\label{equation_Wigner_integral_general}
\end{equation}

Our results for the Cauchy transform of $X_{t}$ can be extended to this more
general setting at the expense of more cumbersome notation. 

The existence of the solution of equation (\ref%
{equation_Wigner_integral_general}) may fail for large $t$ if the norm of the
solution approaches infinity in finite time. However, for sufficiently small 
$t>0$ we have the following local existence result. (See Theorem 3.1 in \cite%
{biane_speicher01} for a sufficient condition of the global existence in a
simpler class of free SDEs, and Theorem 5.2.1 in \cite{oksendal03} for an
existence and uniqueness result in the case of classical SDEs).

\begin{theorem}
Suppose that $a_{i}$, $b_{i},$ and $c_{i}$ are locally operator Lipshitz
functions and $\overline{X}$ is bounded in operator norm. Then, there exist $%
t_{0}>0$ and a family of operators $X_{t}$ defined for all $t\in \left[
0,t_{0}\right) $ and bounded in operator norm, such that $X_{0}=\overline{X}%
, $ and $X_{t}$ is a unique solution of (\ref%
{equation_Wigner_integral_general}) for $t<t_{0}.$
\end{theorem}

\textbf{Proof:} The proof proceeds by Picard's method of successive
approximations. We will give the prove for the case $m=1.$ The general case
is similar. Define $X_{t}^{\left( 0\right) }=\overline{X},$ and 
\begin{equation}
X_{t}^{\left( N+1\right) }=\overline{X}+\int_{0}^{t}a(X_{\tau }^{\left(
N\right) })d\tau +\int_{0}^{t}b(X_{\tau }^{\left( N\right) })dW_{\tau
}c(X_{\tau }^{\left( N\right) }).  \label{Picard_approximation_scheme}
\end{equation}%
We aim to show that this process converges for all sufficiently small $t$.
For this, it is enough to show that for a sufficiently small $t_{0}>0$ $\ $%
and all $t<t_{0}$ and $N\geq 1,$ the following two claims hold: (i) 
\begin{equation*}
\left\| X_{t}^{\left( N\right) }-X_{t}^{\left( N-1\right) }\right\| ^{2}\leq
C\frac{R^{N}}{N!}t^{N},
\end{equation*}%
for some constant $C$ and $R,$ and (ii) 
\begin{equation*}
\left\| X_{t}^{\left( N\right) }\right\| \leq A
\end{equation*}%
for a constant $A.$

Claim (ii) follows from claim (i) because (i) implies that 
\begin{equation*}
\left\| X_{t}^{\left( N\right) }-\overline{X}\right\| \leq Cf\left( t\right)
,
\end{equation*}%
where $f\left( t\right) :=\sum_{k=1}^{\infty }\frac{R^{k/2}}{\sqrt{k!}}%
t^{k/2}$ is defned for all $t<1/R,$ monotonically increasing, differentiable
at $0,$ and vanishes at zero. This implies that for all $A>\left\| \overline{%
X}\right\| ,$ there exists such $t_{0}>0$ that $\left\| X_{t}^{\left(
N\right) }\right\| \leq A$ for all $t<t_{0}.$ Moreover, this choice of $t_{0}
$ is independent of $N.$

In order to prove (i), we proceed by induction. The case $N=1$ is special
and can be easily verified separately. Assume that (i) and (ii) hold for $%
X_{t}^{\left( N\right) }$ and $X_{t}^{(N-1)}$ and let us prove that (i)
holds for $X_{t}^{\left( N+1\right) }.$ We write 
\begin{eqnarray}
\left\| X_{t}^{\left( N+1\right) }-X_{t}^{\left( N\right) }\right\| &\leq
&\left\| \int_{0}^{t}\left[ a(X_{\tau }^{\left( N\right) })-a(X_{\tau
}^{\left( N-1\right) })\right] d\tau \right\|
\label{inequality_Picard_method} \\
&&+\left\| \int_{0}^{t}b(X_{\tau }^{\left( N\right) })dW_{\tau }c(X_{\tau
}^{\left( N\right) })-\int_{0}^{t}b(X_{\tau }^{\left( N-1\right) })dW_{\tau
}c(X_{\tau }^{\left( N-1\right) })\right\| .  \notag
\end{eqnarray}%
The second term in this expression can be estimated by the following sum: 
\begin{eqnarray*}
&&\left\| \int_{0}^{t}\left[ b(X_{\tau }^{\left( N\right) })-b\left( X_{\tau
}^{\left( N-1\right) }\right) \right] dW_{\tau }c(X_{\tau }^{\left( N\right)
})\right\| \\
&&+\left\| \int_{0}^{t}b(X_{\tau }^{\left( N-1\right) })dW_{\tau }\left[
c(X_{\tau }^{\left( N\right) })-c(X_{\tau }^{\left( N-1\right) })\right]
\right\| \\
&\leq &2\sqrt{2}\left( \int_{0}^{t}\left\| b(X_{\tau }^{\left( N\right)
})-b\left( X_{\tau }^{\left( N-1\right) }\right) \right\| ^{2}\left\|
c(X_{\tau }^{\left( N\right) })\right\| ^{2}d\tau \right) ^{1/2} \\
&&+2\sqrt{2}\left( \int_{0}^{t}\left\| c(X_{\tau }^{\left( N\right)
})-c\left( X_{\tau }^{\left( N-1\right) }\right) \right\| ^{2}\left\|
b(X_{\tau }^{\left( N-1\right) })\right\| ^{2}d\tau \right) ^{1/2},
\end{eqnarray*}%
where we used the free Burkholder-Gundy inequality (see Theorem 3.2.1 in %
\cite{biane_speicher98}).

By using the assumption that $b$ and $c$ are operator Lipschitz and claim
(i), we see that this expression is bounded by 
\begin{equation*}
4\sqrt{2}\left[ K_{A}^{2}C\frac{R^{N}}{N!}A^{2}\int_{0}^{t}\tau ^{N}d\tau %
\right] ^{1/2}=\left[ 32K_{A}^{2}A^{2}C\frac{R^{N}}{\left( N+1\right) !}%
t^{N+1}\right] ^{1/2}.
\end{equation*}%
A similar estimate can be obtained for the first part of (\ref%
{inequality_Picard_method}), and by worsening a constant, we can obtain the
following inequality: 
\begin{equation*}
\left\| X_{t}^{\left( N+1\right) }-X_{t}^{\left( N\right) }\right\| ^{2}\leq
64K_{A}^{2}A^{2}C\frac{R^{N}}{\left( N+1\right) !}t^{N+1}
\end{equation*}%
provided that $t<t_{0}.$ This shows that claim (i) holds for $X_{t}^{\left(
N+1\right) }$ provided that $R>64K_{A}^{2}A^{2}.$

Hence, the sequence $X_{t}^{\left( N\right) }=\left( X_{t}^{\left( N\right)
}-X_{t}^{\left( N-1\right) }\right) +\ldots +\left( X_{t}^{\left( 1\right) }-%
\overline{X}\right) +\overline{X}$ is convergent in operator norm for every $%
t<t_{0}$. Let the limit be denoted by $X_{t}.$ By using the free
Burkholder-Gundy inequality, we can take limits on both sides of (\ref%
{Picard_approximation_scheme}) and check that $X_{t}$ is a solution of (\ref%
{equation_Wigner_integral}).

Next, suppose that $X_{t}$ and $X_{t}^{\prime }$ are two different solutions
of (\ref{equation_Wigner_integral}) for $t<t_{0}$. Let $v\left( t\right)
=\left\| X_{t}-X_{t}^{\prime }\right\| .$ By using the assumption that the
coefficients are operator Lipschitz and by using the free Burkholder-Gundy
inequality, we obtain: 
\begin{eqnarray*}
v\left( t\right) &\leq &c_{1}\int_{0}^{t}v\left( \tau \right) d\tau
+c_{2}\left( \int_{0}^{t}v\left( \tau \right) ^{2}d\tau \right) ^{1/2} \\
&\leq &\left( c_{1}\sqrt{t_{0}}+c_{2}\right) \left( \int_{0}^{t}v\left( \tau
\right) ^{2}d\tau \right) ^{1/2},
\end{eqnarray*}%
where $c_{1}$ and $c_{2}$ are certain positive constants that depend on
Lipschitz constants. (The second inequality follows by the Cauchy-Schwarz
inequality.) By the Grownwall inequality (see \cite{oksendal03}, exercise
5.17 on p.80) it follows that $v\left( t\right) ^{2}=0$ for all $t<t_{0}.$
Hence $X_{t}=X_{t}^{\prime }$ and we established the uniqueness of the
solution. QED.

\subsection{Equations for the Cauchy Transform\label{subsection_theory}}

\begin{theorem}
\label{theorem_g_evolution_Wigner} Assume that $a,$ $b,$ and $c$ are locally
operator-Lipschitz functions and let $X_{t}$ be a solution of equation (%
\ref{equation_Wigner}) bounded in operator norm for all $t \in [0,t_0)$. Let $G_{t}$ and $g_{t}$ denote the resolvent of $%
X_{t}$ and the expectation of the resolvent, respectively, and let $%
a_{t}=a(X_{t}),$ $b_{t}=b(X_{t}),$ and $c_{t}=c(X_{t}).$ Then, for all $t \in [0,t_0)$,
\begin{equation}
\frac{dg_{t}}{dt}%
=-E(a_{t}G_{t}^{2})+E(b_{t}c_{t}G_{t})E(b_{t}c_{t}G_{t}^{2}).
\label{equation_main}
\end{equation}
\end{theorem}

Let us mention an important particular case, when the product\ $bc$ does not
depend on $X_{t}.$ In this case, the equation simplifies to the following: 
\begin{equation}
\frac{dg_{t}}{dt}=-E(a_{t}G_{t}^{2})+[E(bc)]^{2}g_{t}\frac{\partial g_{t}}{%
\partial z}.  \label{equation_main_independent_bc}
\end{equation}
If we assume in addition that $b=c=1$ and $a(x)$ is a polynomial then (\ref%
{equation_main_independent_bc}) implies the free Fokker-Planck equation (\ref%
{PDE_free_Fokker_Planck}) of Biane and Speicher. We will demonstrate this in
Corollary \ref{corollary2} below.

In the proof of Theorem \ref{theorem_g_evolution_Wigner} we need the
following lemma.

\begin{lemma}
\label{lemma_Ito_isometry_variant} Let operators $H_{1}$ and $H_{2}$ belong
to the subalgebra $\mathcal{W}_{a}$ which is generated by $\{W_{\tau }\}$
where $\tau \leq a.$ Then%
\begin{eqnarray*}
&&E\left[ \left( \int_{a}^{b}b_{\tau }(dW_{\tau })c_{\tau }\right)
H_{1}\left( \int_{a}^{b}b_{\tau }(dW_{\tau })c_{\tau }\right) H_{2}\right] 
\\
&=&\int_{a}^{b}E(c_{\tau }H_{1}b_{\tau })E(c_{\tau }H_{2}b_{\tau })d\tau .
\end{eqnarray*}
\end{lemma}

This result follows if we write the integral as the limit of sums and use
formula (\ref{formula_4_product})$.$

\textbf{Proof of Theorem}: For conciseness of the following formulas, let us
use the following notation: 
\begin{equation*}
A=\int_{t}^{t+\Delta t}a_{\tau }d\tau ,
\end{equation*}%
and 
\begin{equation*}
B=\int_{t}^{t+\Delta t}b_{\tau }(dW_{\tau })c_{\tau }.
\end{equation*}%
Note that $\Vert A\Vert _{2}=O(\Delta t)$ and $\Vert B\Vert _{2}=O(\sqrt{%
\Delta t})$ for small $\Delta t.$

By using the resolvent identity twice, we can write:%
\begin{eqnarray*}
G_{t+\Delta t}-G_{t} &=&-G_{t+\Delta t}(A+B)G_{t} \\
&=&-G_{t}AG_{t}-G_{t}BG_{t}+G_{t+\Delta t}(A+B)G_{t}(A+B)G_{t}
\end{eqnarray*}%
Note that 
\begin{equation*}
\Vert G_{t+\Delta t}AG_{t}AG_{t}+G_{t+\Delta t}AG_{t}BG_{t}+G_{t+\Delta
t}BG_{t}AG_{t}\Vert _{2}=o(\Delta t).
\end{equation*}%
In addition, 
\begin{equation*}
\Vert G_{t+\Delta t}-G_{t}\Vert _{2}=O(\sqrt{\Delta t}),
\end{equation*}%
which implies 
\begin{equation*}
\left\| G_{t+\Delta t}BG_{t}BG_{t}-G_{t}BG_{t}BG_{t}\right\| _{2}=o\left(
\Delta t\right) 
\end{equation*}%
Hence, we can write 
\begin{equation}
E(G_{t+\Delta
t}-G_{t})=E(-G_{t}AG_{t}-G_{t}BG_{t}+G_{t}BG_{t}BG_{t})+o(\Delta t)
\label{formula_theorem1_deriviation}
\end{equation}%
Next, we use the facts that 
\begin{equation*}
\int_{t}^{t+\Delta t}a_{\tau }d\tau =a_{t}\Delta t+o(\Delta t),
\end{equation*}%
that 
\begin{equation*}
E[\int_{t}^{t+\Delta t}G_{t}b_{\tau }(dW_{\tau })c_{\tau }G_{t}]=0,
\end{equation*}%
and that 
\begin{eqnarray*}
&&E\left[ G_{t}\left( \int_{t}^{t+\Delta t}b_{\tau }(dW_{\tau })c_{\tau
}\right) G_{t}\left( \int_{t}^{t+\Delta t}b_{\tau }(dW_{\tau })c_{\tau
}\right) G_{t}\right]  \\
&=&(\Delta t)E(c_{t}G_{t}b_{t})E(c_{t}G_{t}^{2}b_{t})+o(\Delta t),
\end{eqnarray*}%
where the latter holds because of Lemma \ref{lemma_Ito_isometry_variant} and
the assumption that $b_{t}$ and $c_{t}$ are Lipschitz. Hence, after taking
the expectation in (\ref{formula_theorem1_deriviation}) we obtain 
\begin{equation*}
g_{t+\Delta t}-g_{t}=\Delta
t\{-E(a_{t}G_{t}^{2})+E(c_{t}G_{t}b_{t})E(c_{t}G_{t}^{2}b_{t})\}+o(\Delta t),
\end{equation*}%
which is equivalent to the statement of the theorem. QED.

In order to proceed further and obtain a differential equation on $g_{t}$,
we need to impose additional conditions on $a_{t}$, $b_{t},$ and $c_{t}$
which would allow us to eliminate expectations from (\ref{equation_main})$.$

\begin{proposition}
\label{proposition_g_evolution_PDE}Let $X_{t}$ be the solution of equation (%
\ref{equation_Wigner}), and $G( t,z) $ and $g( t,z) $ be its resolvent and
the expectation of the resolvent, respectively. Suppose that functions $a$
and $bc$ are polynomials in one variable and that their degrees are not
greater than $k\geq 0.$ Then, 
\begin{eqnarray}
\frac{dg}{dt} &=&-\frac{\partial ( ag) }{\partial z}-\sum_{j=0}^{k-2}\frac{(
k-1-j) !}{k!}\frac{\partial ^{j+2}a( z) }{( \partial z) ^{j+2}}E( X^{j}) \\
&&+[ bcg+\sum_{j=0}^{k-1}\frac{( k-1-j) !}{k!}\frac{\partial ^{j+1}[ b( z)
c( z) ] }{( \partial z) ^{j+1}}E( X_{t}^{j}) ] \\
&&\times [ \frac{\partial ( bcg) }{\partial z}+\sum_{j=0}^{k-2}\frac{(
k-1-j) !}{k!}\frac{\partial ^{j+2}[ b( z) c( z) ] }{( \partial z) ^{j+2}}E(
X_{t}^{j}) ] .
\end{eqnarray}
\end{proposition}

This equation is more useful than it might seem at the first sight. First of
all, it is often possible to compute the expectations $E( X_{t}^{j}) $by
using the It\^{o} formula. Second, if these expectations are known, then the
equation is a quasilinear PDE and the method of characteristics is
applicable.

\textbf{Proof:} Let $f(x)$ be a polynomial. If we expand $f(X)(X-z)^{-1}$
and $f(X)(X-z)^{-2}$ in partial fractions and then take the expectations, we
obtain the formulas: 
\begin{equation*}
E(\frac{f(X)}{X-z})=E(\frac{f(z)}{X-z})+\sum_{j=0}^{k-1}\frac{(k-1-j)!}{k!}%
\frac{\partial ^{j+1}f(z)}{(\partial z)^{j+1}}E(X^{j}),
\end{equation*}%
and 
\begin{eqnarray*}
E(\frac{f(X)}{(X-z)^{2}}) &=&E(\frac{f(z)}{(X-z)^{2}})+E(\frac{f^{\prime }(z)%
}{X-z}) \\
&&+\sum_{j=0}^{k-2}\frac{(k-1-j)!}{k!}\frac{\partial ^{j+2}f(z)}{(\partial
z)^{j+2}}E(X^{j}).
\end{eqnarray*}%
By using $a(z)$ or $b(z)c(z)$ as $f(z)$ it is easy to see that the statement
of the proposition follows from Theorem \ref{theorem_g_evolution_Wigner}.
QED.

\begin{corollary}
\label{corollary1}Suppose that $a$ is a polynomial and that $bc=1.$ Then, 
\begin{eqnarray}
\frac{dg}{dt} &=&-\frac{\partial ( ag) }{\partial z}+g\frac{\partial g}{%
\partial z} \\
&&-\sum_{j=0}^{k-2}\frac{( k-1-j) !}{k!}\frac{\partial ^{j+2}a( z) }{(
\partial z) ^{j+2}}E( X^{j})
\end{eqnarray}
\end{corollary}

\bigskip 

\begin{corollary}
\label{corollary2}Suppose that $a$ is a polynomial with real coefficients,
that $b=c=1,$ and that $X_{0}$ is self-adjoint. Assume that the spectral
distribution of $X_{t}$ is absolutely continuous and bounded with the
density $p\left( x,t\right) .$ Then, at all points where $\partial
p/\partial x$ is defined, it is true that%
\begin{equation*}
\frac{dp}{dt}=-\frac{\partial }{\partial x}\{ap+p\cdot Hp\},
\end{equation*}%
where $Hp$ is the Hilbert transform of $p.$
\end{corollary}

(This is the free Fokker-Planck equation (\ref{PDE_free_Fokker_Planck}) of
Biane and Speicher.)

\textbf{Proof of Corollary \ref{corollary2}:} Note that $X_{t}$ are self-adjoint for all $t.$ Let us take
the imaginary part on both sides of the formula in Corollary \ref{corollary1}%
, and then pass to the limit $y\rightarrow 0,$ where $y:=\mathfrak{\func{Im}}%
z$. Assume that $g(z)$ is analytic at $z=x$ and therefore taking the limit
commutes with operations of differentiation with respect to $t$ and $z$.

Since $X_{t}$ is self-adjoint, therefore $\mathfrak{\func{Im}}%
(E(X_{t}^{j}))=0.$ Hence, the formula in Corollary \ref{corollary1}
simplifies as follows:%
\begin{equation}
\pi \frac{dp}{dt}=\lim_{y\rightarrow 0}\{-\frac{\partial }{\partial z}%
\mathfrak{\func{Im}}[ag]+\mathfrak{\func{Im}}[g\frac{\partial g}{\partial z}%
]\},  \label{equation_corollary2_1}
\end{equation}%
where we used the Stieltjes inversion formula. Note that 
\begin{eqnarray*}
\lim_{y\rightarrow 0}\mathfrak{\func{Im}}(g\frac{\partial g}{\partial z})
&=&\lim_{y\rightarrow 0}\{\mathfrak{\func{Im}}g\mathfrak{\func{Re}}\frac{%
\partial g}{\partial z}+\mathfrak{\func{Re}}g\mathfrak{\func{Im}}\frac{%
\partial g}{\partial z}\} \\
&=&\pi (p\frac{\partial }{\partial x}(-Hp)-\left( Hp\right) \frac{\partial }{%
\partial x}p) \\
&=&-\pi \frac{\partial }{\partial x}\left[ p\cdot Hp\right] .
\end{eqnarray*}%
Similarly, 
\begin{equation*}
\lim_{y\rightarrow 0}\frac{\partial }{\partial z}\mathfrak{\func{Im}}%
[ag]=\pi \frac{\partial }{\partial x}\left[ ap\right] 
\end{equation*}%
because $\mathfrak{\func{Im}}a(x)=0$ and $\mathfrak{\func{Re}}a(x)=a(x).$
Hence, equation (\ref{equation_corollary2_1}) simplifies to 
\begin{equation*}
\frac{dp}{dt}=-\frac{\partial }{\partial x}\left[ ap+p\cdot Hp\right] .
\end{equation*}

QED.

\subsection{Examples\label{subsection_examples}}

In this section, we calculate explicit solutions in several particular cases.

\subsubsection{Ornstein-Uhlenbeck}

\begin{proposition}
Suppose that $X_{t}$ satisfies the equation of the free Ornstein-Uhlenbeck
process:%
\begin{equation*}
dX_{t}=\theta X_{t}dt+\sigma dW_{t}.
\end{equation*}%
Suppose that $X_{0}=0$. Then the spectral distribution of $X_{t}$ is the
semicircle distribution supported on the interval $I_{\theta },$ where $%
\newline
$1) 
\begin{equation*}
I_{\theta }=[ -\sqrt{\frac{2\sigma ^{2}}{\theta }( e^{2\theta t}-1) },+\sqrt{%
\frac{2\sigma ^{2}}{\theta }( e^{2\theta t}-1) }] ,
\end{equation*}%
if $\theta >0$, \newline
2) 
\begin{equation*}
I_{\theta }=[ -2\sigma \sqrt{t},2\sigma \sqrt{t}]
\end{equation*}%
if $\theta =0,$ and\newline
3) 
\begin{equation*}
I_{\theta }=[ -\sqrt{\frac{2\sigma ^{2}}{| \theta | }( 1-e^{-2| \theta | t}) 
},+\sqrt{\frac{2\sigma ^{2}}{| \theta | }( 1-e^{-2| \theta | t}) }] ,
\end{equation*}%
if $\theta <0.$
\end{proposition}

Hence, if $\theta >0,$ then the support of the distribution grows
exponentially; if $\theta =0,$ then the support grows linearly, and if $%
\theta <0,$ the spectral distribution converges to the semicircle
distribution supported on the interval $[ -\sigma \sqrt{2/\vert \theta \vert 
},\sigma \sqrt{2/\vert \theta \vert }] .$

\textbf{Proof}: In this case $a_{t}=\theta X_{t},$ $b_{t}=\sqrt{\sigma }.$
Note that 
\begin{equation*}
a_{t}=\theta X_{t}=\theta (z+G_{t}^{-1}).
\end{equation*}%
Hence, 
\begin{equation*}
a_{t}G_{t}^{2}=\theta (zG_{t}^{2}+G_{t})
\end{equation*}%
and 
\begin{equation*}
E(a_{t}G_{t}^{2})=\theta (z_{t}\frac{\partial g_{t}}{\partial z}+g_{t}).
\end{equation*}%
Therefore, the differential equation for $g_{t}$ is 
\begin{equation}
\frac{\partial g}{\partial t}+(\theta z-\sigma ^{2}g)\frac{\partial g}{%
\partial z}=-\theta g_{t}.  \label{equation_g_W}
\end{equation}%
The initial condition $X_{0}=0$ corresponds to $g(0,z)=-z^{-1},$ and we can
solve this partial differential equation by using the method of
characteristics (see pp. 9-19 in \cite{john81}).

Indeed the equations of characteristic curves are 
\begin{eqnarray}
\frac{dt}{d\xi } &=&1,  \label{equations_case1_1} \\
\frac{dz}{d\xi } &=&\theta z-\sigma ^{2}g,  \label{equations_case1_2} \\
\frac{dg}{d\xi } &=&-\theta g.  \label{equations_case1_3}
\end{eqnarray}

By using (\ref{equations_case1_1}), we can set $\xi =t.$ Then (\ref%
{equations_case1_3}) implies that 
\begin{equation*}
g(t)=Ae^{-\theta t},
\end{equation*}%
and then we can solve (\ref{equations_case1_2}) as 
\begin{equation*}
z(t)=Ce^{\theta t}+\frac{\sigma ^{2}A}{2\theta }e^{-\theta t}.
\end{equation*}%
It follows that the initial point of a characteristic curve is given by
equations: 
\begin{eqnarray*}
g(0) &=&A, \\
z(0) &=&C+\frac{\sigma ^{2}A}{2\theta }.
\end{eqnarray*}%
On the other hand we can parameterize the initial condition of the PDE as
follows: 
\begin{equation*}
z(s)=s,\;g(s)=-1/s.
\end{equation*}%
Hence, we obtain the following paramterization for $A$ and $C$:%
\begin{equation*}
A=-1/s,\;C=s+\frac{\sigma ^{2}}{2\theta }\frac{1}{s}.
\end{equation*}%
Therefore, the equations of the characteristic surface are 
\begin{eqnarray}
g(s,t) &=&-\frac{1}{s}e^{-\theta t},  \label{equations_case1_4} \\
z(s,t) &=&(s+\frac{\sigma ^{2}}{2\theta }\frac{1}{s})e^{\theta t}-\frac{%
\sigma ^{2}}{2\theta }\frac{1}{s}e^{-\theta t}.  \label{equations_case1_5}
\end{eqnarray}%
From (\ref{equations_case1_4}) we have 
\begin{equation*}
s=-\frac{1}{ge^{\theta t}}.
\end{equation*}%
After we substitute this in (\ref{equations_case1_5}) and re-arrange, we
obtain the following functional equation for $g(t,z)$:%
\begin{equation*}
g^{2}+\frac{2\theta z}{\sigma ^{2}(e^{2\theta t}-1)}g+\frac{2\theta }{\sigma
^{2}(e^{2\theta t}-1)}=0,
\end{equation*}%
provided that $\theta \neq 0.$ We can easily solve this quadratic equation
for $g.$ Note that by the Stieltjes inversion formula the density of the
corresponding distribution is given by the imaginary part of the Cauchy
transform $g.$ We can check that in our case this density corresponds to the
density of the semicircle distribution. The radius of the semicircle
distribution is 
\begin{equation*}
\sqrt{\frac{2\sigma ^{2}}{\theta }(e^{2\theta t}-1)},
\end{equation*}%
if $\theta >0,$ and 
\begin{equation*}
\sqrt{\frac{2\sigma ^{2}}{|\theta |}(1-e^{-2|\theta |t})},
\end{equation*}%
if $\theta <0.$ This implies the statement of the proposition for $\theta
\neq 0.$ The case $\theta =0$ can be analyzed similarly. QED.

\subsubsection{Geometric Brownian Motion}

Now let us consider the case when the coefficient $b_{t}$ explicitly depends
on $X_{t}.$ Namely, let $a_{t}=\theta X_{t},$ and $b_{t}=X_{t}^{1/2}.$

In this example we deal with the equation 
\begin{equation*}
dX_{t}=\theta X_{t}dt+X_{t}^{1/2}( dW_{t}) X_{t}^{1/2},
\end{equation*}%
which is an analog of the classical equation for the ``geometric Brownian
motion'', $dx_{t}=\theta x_{t}dt+x_{t}dB_{t}.$

Let us assume that $X_{0}=I$ and use the free Ito formula to study the
moments of the solution. Clearly, $E( X_{t}) =e^{\theta t}.$ In order to
calculate the second moment, we write 
\begin{eqnarray*}
d( X_{t}^{2}) &=&( X_{t}+dX_{t}) ^{2}-X_{t}^{2} \\
&=&( 2\theta X_{t}^{2}+e^{\theta t}X_{t}) dt+X_{t}^{3/2}( dW_{t})
X_{t}^{1/2}+X_{t}^{1/2}( dW_{t}) X_{t}^{3/2},
\end{eqnarray*}%
where we used the free Ito formula to calculate 
\begin{eqnarray*}
dX_{t}dX_{t} &=&X_{t}^{1/2}( dW_{t}) X_{t}( dW_{t}) X_{t}^{1/2} \\
&=&E( X_{t}) X_{t}dt
\end{eqnarray*}

Let $h_{t}$ denote $E( X_{t}^{2}) .$ Then we have the following equation: 
\begin{equation*}
\frac{dh_{t}}{dt}=2\theta h_{t}+e^{2\theta t}
\end{equation*}%
with the initial condition $h_{0}=1.$ The solution is 
\begin{equation*}
h_{t}=( t+1) e^{2\theta t}.
\end{equation*}

Hence, the variance of the spectral distribution of $X_{t}$ is $te^{2\theta
t}.$ The ratio of the standard deviation to the expectation of $X_{t}$ is $%
\sqrt{t}.$

In order to recover the entire spectral distribution, we use Theorem \ref%
{theorem_g_evolution_Wigner}, and obtain the following result.

\begin{proposition}
\label{proposition_example2}Suppose that $X_{t}$ satisfies the following
equation:%
\begin{equation*}
dX_{t}=\theta X_{t}dt+X_{t}^{1/2}(dW_{t})X_{t}^{1/2},
\end{equation*}%
and that $X_{0}=I$. Then, the expectation of the resolvent satisfies the
following functional equation:%
\begin{equation}
z+g^{-1}=e^{(\alpha -zg)t},  \label{equation_geom_Brown}
\end{equation}%
where $\alpha =\theta -1.$ The density of the spectral distribution of $%
X_{t} $ is supported on the interval%
\begin{equation*}
I=[\frac{r_{1}(t)}{1+r_{1}(t)}e^{(\alpha -r_{1}(t))t},\frac{r_{2}(t)}{%
1+r_{2}(t)}e^{(\alpha -r_{2}(t))t}],
\end{equation*}%
where 
\begin{equation*}
r_{1,2}(t)=\frac{-1\pm \sqrt{1+4/t}}{2}.
\end{equation*}
\end{proposition}

We can see from this proposition that the solution of the free SDE exists
remains positive definite for all $t>0$.

If $t\rightarrow \infty ,$ then $r_{1,2}( t) $ are asymptotically $1/t$ and $%
-1-1/t.$ Hence, as $t\rightarrow \infty $ the support of the solution
becomes asymptotically close to 
\begin{equation*}
[ \frac{1}{et}e^{( \theta -1) t},ete^{\theta t}] .
\end{equation*}

In particular, if $\theta <0,$ then both the lower and the upper bound of
the spectral distribution shrink to zero exponentially fast, although at
different rates ($\theta -1$ and $\theta $). If $\theta =0,$ then the lower
bound shrinks to zero exponentially and the upper bound grows linearly. If $%
\theta \in ( 0,1) ,$ then the lower bound shrinks exponentially and the
upper bound grows exponentially. If $\theta =1,$ then the lower bound
declines as $( et) ^{-1}$ and the upper bound grows exponentially. If $%
\theta >1,$ then both the upper and lower bounds grow exponentially.

\textbf{Proof of Proposition \ref{proposition_example2}}: We have 
\begin{equation*}
E( b_{t}^{\ast }G_{t}b_{t}) =E( G_{t}X_{t}) =1+zg_{t},
\end{equation*}%
and 
\begin{equation*}
E( b_{t}^{\ast }G_{t}^{2}b_{t}) =E( G_{t}^{2}X_{t}) =g_{t}+z\frac{\partial
g_{t}}{\partial z}.
\end{equation*}%
Hence the differential equation is 
\begin{equation*}
\frac{\partial g}{\partial t}=-\theta ( g+z\frac{\partial g}{\partial z}) +(
1+zg_{t}) ( g_{t}+z\frac{\partial g_{t}}{\partial z}) ,
\end{equation*}%
or 
\begin{equation}
\frac{\partial g}{\partial t}+z( ( \theta -1) -zg) \frac{\partial g}{%
\partial z}=-g( ( \theta -1) -zg) .  \label{equations_case2_0}
\end{equation}

By assumption, the initial condition is $g( 0,z) =( 1-z) ^{-1}$. \ 

The equations of characteristic curves are 
\begin{eqnarray}
\frac{dt}{d\xi } &=&1,  \label{equations_case2_1} \\
\frac{dz}{d\xi } &=&z( \theta -1-zg) ,  \label{equations_case2_2} \\
\frac{dg}{d\xi } &=&-g( \theta -1-zg) .  \label{equations_case2_3}
\end{eqnarray}
From (\ref{equations_case2_1}) we can set $\xi =t.$ Then, if we divide (\ref%
{equations_case2_3}) by (\ref{equations_case2_2}), we obtain the following
equation: 
\begin{equation*}
\frac{dg}{dz}=-\frac{g}{z},
\end{equation*}%
which implies the following family of equations for the characteristic
curves. 
\begin{equation*}
g=\frac{A}{z}.
\end{equation*}

If we substitute this in equation (\ref{equations_case2_2}) and solve the
resulting ODE, we find:%
\begin{equation*}
z( t) =Ce^{( \theta -1-A) t}.
\end{equation*}%
Hence, 
\begin{equation*}
g( t) =(A/C)e^{-( \theta -1-A) t}.
\end{equation*}

In particular, $z( 0) =C,$ $g( 0) =A/C.$

On the other hand, we can parameterize the initial condition of (\ref%
{equations_case2_0}) as follows:%
\begin{equation*}
z(s)=s,\;g( s) =\frac{1}{1-s}.
\end{equation*}

This implies the following parameterization for $A$ and $C$: 
\begin{equation*}
C=s,\;A=\frac{s}{1-s}.
\end{equation*}%
Hence, the characteristic surface is 
\begin{eqnarray}
z( s,t) &=&s\exp \{ ( \theta -1-\frac{s}{1-s}) t\} ,
\label{equations_case2_4} \\
g( s,t) &=&\frac{1}{1-s}\exp \{ -( \theta -1-\frac{s}{1-s}) t\} .
\label{equations_case2_5}
\end{eqnarray}%
We can eliminate $s$ from these equations: 
\begin{equation*}
s=\frac{zg}{1+zg}.
\end{equation*}%
After we substitute this expression for $s$ in (\ref{equations_case2_4}) and
re-arrange the terms, then we obtain the following equation:%
\begin{equation*}
z+g^{-1}=\exp \{ ( \theta -1-zg) t\} .
\end{equation*}

Let us denote $\theta -1$ as $\alpha $ for simplicity of notation. Then, the
functional equation for the Cauchy transform $g( t,z) $ is as follows: 
\begin{equation}
z+g^{-1}=e^{( \alpha -zg) t}.  \label{formula_example2}
\end{equation}

If we take the differential of this equation, then we find that 
\begin{equation*}
dz(1+gte^{(\alpha -zg)t})=dg(g^{-2}-zte^{(\alpha -zg)t}).
\end{equation*}%
The branch points of the function $g(z)$ can be found from the equation $%
dz/dg=0.$ Hence, at the branch points,%
\begin{equation*}
e^{(\alpha -zg)t}=\frac{1}{g^{2}zt}.
\end{equation*}%
Substituting this into equation (\ref{formula_example2}), we obtain the
following equation for the branch points:%
\begin{equation*}
t(zg)^{2}+t(zg)-1=0.
\end{equation*}%
Hence, 
\begin{equation}
zg=\frac{-1\pm \sqrt{1+4/t}}{2}\equiv r_{1,2}(t).  \label{equation_r12}
\end{equation}%
Then (\ref{formula_example2}) and (\ref{equation_r12}) imply that at the
branch points, 
\begin{equation*}
g_{1,2}=(1+r_{1,2}(t))e^{-(\alpha -r_{1,2}(t))t}
\end{equation*}%
and 
\begin{equation*}
z_{1,2}=\frac{r_{1,2}(t)}{1+r_{1,2}(t)}e^{(\alpha -r_{1,2}(t))t}.
\end{equation*}%
Finally, note that branch points of the Cauchy transform are bounds for the
support of the spectral probability distribution. QED.

\subsubsection{Geometric Brownian Motion II}

In our next example, we consider a different analog of the classical
geometric Brownian motion equation, namely, the following free SDE:%
\begin{equation*}
dX_{t}=\theta X_{t}dt+X_{t}dW_{t}+(dW_{t})X_{t}.
\end{equation*}%
As in the previous example, assume that $X_{0}=I,$ and note that $%
E(X_{t})=e^{\theta t},$ the same as in the previous example.

It is possible to write a PDE for the expectation of the resolvent in this
example similar to equations in Theorem \ref{theorem_g_evolution_Wigner} and
Proposition \ref{proposition_g_evolution_PDE}. However, it seems that it is
difficult to find an explicit solution of this equation and recover the
spectral distribution function of $X_{t}.$

Still, it is possible to see that the behavior of the solution is quite
different from the behavior of the solution in the previous example by
studying the variance of the solution. By using the free Ito formula, we can
write:%
\begin{eqnarray*}
d(X_{t}^{2}) &=&[2\theta X_{t}^{2}+X_{t}^{2}+2E(X_{t})X_{t}+E(X_{t}^{2})]dt
\\
&&+X_{t}^{2}dW_{t}+2X_{t}(dW_{t})X_{t}+(dW_{t})X_{t}^{2}.
\end{eqnarray*}%
Let $h_{t}$ denote $E(X_{t}^{2}).$ Then we have the following ODE for $h_{t}$%
:%
\begin{equation*}
\frac{dh_{t}}{dt}=2(\theta +1)h_{t}+2e^{2\theta t}.
\end{equation*}%
The initial condition is $h_{0}=1$ and the solution is%
\begin{equation*}
h_{t}=2e^{2(\theta +1)t}-e^{2\theta t}.
\end{equation*}%
Hence the variance of $X_{t}$ is $2e^{2\theta t}(e^{2t}-1),$ and the ratio
of the standard deviation to the expectation is $\sqrt{2(e^{2t}-1)}.$ This
ratio grows exponentially fast with $t,$ quite unlike the previous example,
where this ratio equals $\sqrt{t}.$

\subsubsection{Explosive equation}

In our final example, we will consider an equation whose solution explodes
in finite time. By this we mean that the norm of the solution becomes
infinite in finite time.

\begin{proposition}
\label{proposition_explosive}Suppose that $X_{t}$ satisfies the following
equation:%
\begin{equation*}
dX_{t}=kX_{t}(dW_{t})X_{t},
\end{equation*}%
and let the initial condition be $X_{0}=aI.$ Then the spectral distribution
of $X_{t}$ is defined for all $t\leq (ak)^{-2}$ and it is supported on the
interval:%
\begin{equation*}
I=[\frac{(1-ak\sqrt{t})^{2}}{(1-a^{2}k^{2}t)^{2}},\frac{(1+ak\sqrt{t})^{2}}{%
(1-a^{2}k^{2}t)^{2}}].
\end{equation*}%
For $\tau \in (0,1),$ the density of the spectral distribution of the
operator $a^{-1}X_{(ak)^{2}\tau }$ is given by the formula:%
\begin{equation*}
f(\xi )=\frac{\sqrt{-(1-\tau )^{2}\xi ^{2}+2(1+\tau )\xi -1}}{2\pi \xi
^{3}\tau }.
\end{equation*}
\end{proposition}

\textbf{Proof of Proposition \ref{proposition_explosive}:} We can compute 
\begin{eqnarray*}
E(G_{t}b_{t}^{2}) &=&k[E(G_{t}^{-1})+2z+z^{2}g_{t}] \\
&=&k[a+z+z^{2}g_{t}],
\end{eqnarray*}%
where we used the fact that 
\begin{equation*}
E(G_{t}^{-1})=E(X_{t})-z=E(X_{0})-z=a-z.
\end{equation*}%
In addition, 
\begin{equation*}
E(G_{t}^{2}b_{t}^{2})=k[1+2zg_{t}+z^{2}\frac{\partial g_{t}}{\partial z}].
\end{equation*}

Hence, the differential equation for $g_{t}$ is 
\begin{equation*}
\frac{\partial g}{\partial t}-k^{2}( a+z+z^{2}g) z^{2}\frac{\partial g}{%
\partial z}=k^{2}( a+z+z^{2}g) ( 1+2zg) ,
\end{equation*}%
and the initial condition is $g_{0}=( a-z) ^{-1}.$

The equations for the characteristic curves are%
\begin{eqnarray}
\frac{dt}{d\xi } &=&1,  \label{equations_case4_1} \\
\frac{dz}{d\xi } &=&-k^{2}( a+z+z^{2}g) z^{2},  \label{equations_case4_2} \\
\frac{dg}{d\xi } &=&k^{2}( a+z+z^{2}g) ( 1+2zg) .  \label{equations_case4_3}
\end{eqnarray}

From (\ref{equations_case4_1}), we can set $\xi =t,$ and the equations for
the characteristic curves in $( z,g) $-plane become: 
\begin{eqnarray}
\frac{dz}{dt} &=&-k^{2}( a+z+z^{2}g) z^{2},  \label{equations_case4_4} \\
\frac{dg}{dt} &=&k^{2}( a+z+z^{2}g) ( 1+2zg) .  \label{equations_case4_5}
\end{eqnarray}

After dividing (\ref{equations_case4_5}) by (\ref{equations_case4_4}), we
obtain: 
\begin{equation*}
\frac{dg}{dz}=-\frac{1+2gz}{z^{2}}.
\end{equation*}%
The general solution of this equation is 
\begin{equation}
g( z) =-z^{-1}+Cz^{-2}.  \label{equations_case4_6}
\end{equation}%
If we substitute this expression in (\ref{equations_case4_4}), we obtain: 
\begin{equation*}
\frac{dz}{dt}=-k^{2}( a+C) z^{2}.
\end{equation*}%
Hence, 
\begin{equation}
z( t) =\frac{1}{k^{2}( a+C) t+A}.  \label{equations_case4_7}
\end{equation}

By substituting this in (\ref{equations_case4_6}), we obtain: 
\begin{equation}
g( t) =-( a+C) k^{2}t-A+C( ( a+C) k^{2}t+A) ^{2}.  \label{equations_case4_8}
\end{equation}

In particular, if $t=0,$ then%
\begin{eqnarray}
z( 0) &=&1/A,  \label{equations_case4_9} \\
g( 0) &=&-A+CA^{2}.  \notag
\end{eqnarray}

On the other hand, the initial condition is $g( z) =( a-z) ^{-1},$ which we
can parameterize as follows: 
\begin{equation}
z( s) =s,\;g( s) =( a-s) ^{-1}.  \label{equations_case4_10}
\end{equation}

Comparing (\ref{equations_case4_9}) and (\ref{equations_case4_10}), we
obtain the following parameterization for $A$ and $C$:%
\begin{equation*}
A=\frac{1}{s},\;C=\frac{as}{a-s}.
\end{equation*}%
We substitute these expressions in (\ref{equations_case4_7}) and (\ref%
{equations_case4_8}) and obtain:%
\begin{eqnarray}
z(t,s) &=&\frac{1}{\frac{a^{2}}{a-s}k^{2}t+\frac{1}{s}},
\label{equations_case4_11} \\
g(t,s) &=&\frac{a^{5}s}{(a-s)^{3}}k^{4}t^{2}+\frac{a^{2}(a+s)}{(a-s)^{2}}%
k^{2}t+\frac{1}{a-s}.  \label{equations_case4_12}
\end{eqnarray}%
We are going to eliminate $s$ from the pair of equations (\ref%
{equations_case4_11}) and (\ref{equations_case4_12}). For this reason, we
write (\ref{equations_case4_12}) as follows: 
\begin{equation*}
g(t,s)=\frac{s}{a-s}(\frac{a^{2}}{a-s}k^{2}t+\frac{1}{s})(\frac{a^{3}}{a-s}%
k^{2}t+1),
\end{equation*}%
and then we substitute (\ref{equations_case4_11}) and obtain: 
\begin{equation*}
g(t,s)=\frac{s}{a-s}\frac{1}{z}(\frac{a^{3}}{a-s}k^{2}t+1),
\end{equation*}%
or 
\begin{equation}
\frac{a-s}{s}zg=\frac{a^{3}}{a-s}k^{2}t+1.  \label{equations_case4_13}
\end{equation}%
By using (\ref{equations_case4_11}) again, we note that 
\begin{equation*}
\frac{a^{3}}{a-s}k^{2}t+1=\frac{a}{z}-\frac{a}{s}+1.
\end{equation*}%
Hence, (\ref{equations_case4_13}) can be re-written as follows: 
\begin{equation*}
(a-s)zg=(\frac{a}{z}+1)s-a,
\end{equation*}%
and, therefore, 
\begin{equation*}
s=a\frac{1+zg}{1+zg+\frac{a}{z}},
\end{equation*}%
and 
\begin{equation*}
a-s=a\frac{\frac{a}{z}}{1+zg+\frac{a}{z}}.
\end{equation*}%
After substituting these expressions in equation (\ref{equations_case4_11}),
we obtain: 
\begin{equation*}
\frac{1}{z}=(z+z^{2}g+a)k^{2}t+\frac{1}{a}+\frac{1}{z+z^{2}g}.
\end{equation*}%
After re-arranging the terms and dividing by $z,$ we get the following
equation: 
\begin{equation*}
k^{2}tz^{3}g^{2}+(\frac{z}{a}-1+(a+2z)zk^{2}t)g+\frac{1}{a}+(z+a)k^{2}t=0.
\end{equation*}%
This functional equation for $g(z,t)$ is quadratic and therefore it is
easily solvable.

In particular, the branch points of $g_{t}( z) $ are the zeros of the
discriminant of this equation, which can be computed as 
\begin{equation*}
D=( k^{2}at-\frac{1}{a}) ^{2}z^{2}-2( k^{2}at+\frac{1}{a}) z+1.
\end{equation*}%
Therefore, the branch points are 
\begin{equation*}
z_{\pm }=a\frac{( 1\pm ak\sqrt{t}) ^{2}}{( 1-a^{2}k^{2}t) ^{2}}.
\end{equation*}%
Note that as $t$ approaches $( ak) ^{-2},$ the branch points approach $a/4$
and $\infty$.

It follows that for $t<( ak) ^{-2},$ the spectral distribution of $X_{t}$ is
supported on the interval $[ z_{-},z_{+}] $ and in this region it has the
density 
\begin{equation*}
f( x) dx=\frac{1}{a}\frac{\sqrt{-( 1-k^{2}a^{2}t) ^{2}( \frac{x}{a}) ^{2}+2(
k^{2}a^{2}t+1) ( \frac{x}{a}) -1}}{2\pi k^{2}a^{2}t( x/a) ^{3}}dx.
\end{equation*}%
If we use variables $\tau =k^{2}a^{2}t,$ and $\xi =x/a,$ then we can write
this density as%
\begin{equation*}
f( \xi ) d\xi =\frac{\sqrt{-( 1-\tau ) ^{2}\xi ^{2}+2( 1+\tau ) \xi -1}}{%
2\pi \xi ^{3}\tau }.
\end{equation*}%
QED.


\begin{thebibliography}{CD-M}
\bibitem[A]{anshelevich02} Michael Anshelevich. \newblock Ito Formula for
Free Stochastic Integrals. \newblock {\em Journal of Functional Analysis},
188:292--315, 2002.

\bibitem[B]{biane97} Philippe Biane. \newblock Free \mbox{B}rownian motion,
free stochastic calculus and random matrices. \newblock In Dan-Virgil
Voiculescu, editor, \emph{Free Probability Theory}, volume~12 of \emph{%
Fields Institute Communications}, pages 1--19. American Mathematical
Society, 1997.

\bibitem[BSa]{biane_speicher98} Philippe Biane and Roland Speicher. %
\newblock Stochastic calculus with respect to free \mbox{B}rownian motion
and analysis on \mbox{W}igner space. 
\newblock {\em Probability Theory and
Related Fields}, 112:373--409, 1998.

\bibitem[BSb]{biane_speicher01} Philippe Biane and Roland Speicher. \newblock
Free diffusions, free entropy and free \mbox{F}isher information. \newblock  
\emph{Ann. I. H. Poincare}, 37:581--606, 2001.

\bibitem[CD-M]{capitaine_donati-martin05} M.~Capitaine and C.~Donati-Martin. %
\newblock Free \mbox{W}ishart processes. 
\newblock {\em Journal of
Theoretical Probability}, 18:413--438, 2005.

\bibitem[D]{demni08} N.~Demni. \newblock Free \mbox{J}acobi process. %
\newblock {\em Journal of Theoretical Probability}, 21:118--143, 2008.

\bibitem[G]{gao06} Mingchu Gao. \newblock Free \mbox{O}rnstein-\mbox{U}%
hlenbeck processes. 
\newblock {\em Journal of
Mathematical Analysis and Applications}, 322:177--192, 2006.

\bibitem[J]{john81} Fritz John. \newblock Partial Differential Equations. %
\newblock volume 1 of \emph{Applied Mathematical Sciences}. Springer-Verlag,
1981.

\bibitem[KS]{kummerer_speicher92} B.~Kummerer and R.~Speicher. \newblock %
Stochastic integration on the \mbox{C}untz algebra $O_\infty$. \newblock     
\emph{Journal of Functional Analysis}, 103:372--408, 1992.

\bibitem[MP]{marchenko_pastur67} V.A Marcenko and L.A Pastur. \newblock %
Distribution of eigenvalues of some sets of random matrices. 
\newblock {\em
Mathematics in U.S.S.R}, 1:507--536, 1967.

\bibitem[NS]{nica_speicher06} Alexandru Nica and Roland Speicher. \newblock %
Lectures on the combinatorics of free probability. \newblock volume 335 of 
\emph{London Mathematical Society Lecture Note Series}. Cambridge University
Press, 2006.

\bibitem[O]{oksendal03} Bernt Oksendal. \newblock Stochastic Differential
Equations. Sixth Edition. \newblock Springer, 2003.

\bibitem[S]{speicher90} Roland Speicher. \newblock A new example of
independence and white noise. 
\newblock {\em Probability Theory and Related
Fields}, 84:141--159, 1990.

\bibitem[VDN]{voiculescu_dykema_nica92} D.~Voiculescu, K.~Dykema, and
A.~Nica. \newblock {\em Free Random Variables}. \newblock A.M.S. Providence,
RI, 1992. \newblock CRM Monograph series, No.1.
\end{thebibliography}
\end{document}